\documentclass[11pt]{article}
\usepackage[T1]{fontenc}
\usepackage{lmodern,amsmath,amssymb,amsthm}
\usepackage[a4paper,margin=25mm]{geometry}
\usepackage{microtype,enumitem,booktabs}
\usepackage{tikz,float}
\usepackage[hidelinks]{hyperref}
\usepackage{url,authblk,color,cite}	

\setlist[enumerate]{label=\textup{(\roman*)},leftmargin=*}

\newtheorem{theorem}{Theorem}[section]
\newtheorem{lemma}[theorem]{Lemma}
\theoremstyle{remark}
\newtheorem*{remark*}{Remark}
\newcommand{\G}{\Gamma}
\newcommand{\tr}{\operatorname{tr}}
\newcommand{\GDD}{\operatorname{GDDDP}}

\title{$2$-Walk-regular graphs with $k_2$ at most twice the valency}

\author[a]{Zhi Qiao}
\author[b]{Zheng-Jiang Xia 
	{\thanks{Corresponding author. \protect \\
	\indent \ \ 
	E-mail: 
    zhiqiao@sicnu.edu.cn (Z. Qiao),
    xzj@mail.ustc.edu.cn (Z.-J. Xia).}}}
\affil[a]{ School of Mathematical Sciences, Sichuan Normal University, Chengdu, 610068, China}
\affil[b]{ School of Finance, Anhui University of Finance \& Economics, Bengbu, 233030, China}

\date{}

\begin{document} 
\maketitle

\footnote{Dedicated to Professor Jack Koolen on the occasion of his 60th birthday.}

\begin{abstract}
In this paper, we characterize $2$-walk-regular graphs with valency $k$, diameter at least four and $k_2\leq 2k$, generalizing a result of Koolen and Park, and classify the distance-regular graphs satisfying these conditions.
We also classify $2$-walk-regular graphs with valency $k$, diameter at least four and at most $4k+4$ vertices. 
\end{abstract}

\section{Introduction} 
A graph is \emph{$2$-walk-regular} if, for every nonnegative integer
$\ell$, the number of walks of length $\ell$ between vertices at
distance at most two depends only on their distance. Such graphs are generalizations of distance-regular graphs. 
In this paper, we study $2$-walk-regular graphs with valency $k\geq 3$, diameter $D\geq 4$ and $k_2\leq 2k$. 

We first show that such graphs are bipartite if $c_2\geq 2$. 
\begin{theorem}\label{thm:bip}
Let $\G$ be a $2$-walk-regular graph with valency $k\geq3$
and diameter $D\geq4$.
If $k_2\leq 2k$ and $c_2\geq 2$, then $\G$ is bipartite.
\end{theorem}
\begin{remark*}
\begin{enumerate}
    \item Koolen and Park \cite[Theorem 17]{KP12} showed that for $\varepsilon>0$, there exists a real number $\kappa(\varepsilon)$  such that a distance-regular graph with valency $k\geq \kappa(\varepsilon)$, diameter $D\geq 3$ and $k_2\leq (2-\varepsilon)k$ is either a Taylor graph or a bipartite graph with diameter $3$. 
    \item The Hamming graph $H(3,3)$ ($D=3$, $c_2=2$), the line graph of the Petersen graph ($D=3$, $c_2=1$), and the Coxeter graph ($D=4$, $c_2=1$) are nonbipartite distance-regular graphs with $k_2=2k$. 
\end{enumerate}
\end{remark*}

As an application of Theorem \ref{thm:bip}, we obtain the following characterization of $2$-walk-regular graphs with valency $k\geq 3$, diameter $D\geq 4$ and $k_2\leq 2k$. 
\begin{theorem}\label{thm:2k}
Let $\G$ be a $2$-walk-regular graph with valency
$k\geq3$ and diameter $D\geq4$.
Then $k_2\leq2k$ if and only if $\G$ is one of the following graphs:
\begin{enumerate}[ref=\roman*]
\item\label{item:had}
a Hadamard graph, with $k\geq4$ even, $D=4$, and $c_2=k/2$;
\item\label{item:q5}
the $5$-cube, with $k=5$, $D=5$, and $c_2=2$;
\item\label{item:gdd}
the incidence graph of a $\GDD(2,k+1;k;0,\frac{k-1}{2})$,
with $k\geq5$ odd, $D=4$, and $c_2=\frac{k-1}{2}$;
\item\label{item:cubic}
a cubic $1$-walk-regular graph with diameter at least four
and girth at least five; here $k=3$ and $c_2=1$;
\item\label{item:line}
a line graph $L(H)$, where $H$ is a
cubic $3$-walk-regular graph with diameter at least four and girth at least six;
here $k=4$ and $a_1=c_2=1$.
\end{enumerate}
\end{theorem}

\begin{remark*}
\begin{enumerate}
    \item Qiao, Park and Koolen \cite[Theorem 1.3]{QPK19} classified $2$-walk-regular graphs with valency $k\geq 3$, diameter $D\geq 4$ and $c_2>\frac{1}{3}k$. 
    \item Conder and Nedela \cite[Proposition~2.5]{CN07} constructed infinitely many
finite cubic $3$-arc-transitive graphs $H$ with girth at least 10. 
For each such graph $H$, the graphs $H$ and $L(H)$ belong to \eqref{item:cubic} and \eqref{item:line} respectively. 
And Lemma \ref{lem:line} implies that the graphs $L(H)$ are pairwise nonisomorphic. 
Hence both \eqref{item:cubic} and \eqref{item:line} contain infinitely many graphs. 
\end{enumerate}
\end{remark*}

Then we classify distance-regular graphs with valency $k$, diameter $D\geq 4$ and $k_2\leq 2k$. 
\begin{theorem}\label{thm:drg}
Let $\G$ be a distance-regular graph with valency
$k\geq 3$ and diameter $D\geq4$. Then $k_2\leq2k$ if and only if
$\G$ is a Hadamard graph or one of the eleven graphs in
Table~\ref{tab:drg-list}. 
\end{theorem}

\begin{table}[H]
\centering
\caption{Non-Hadamard distance-regular graphs with $k\geq3$, $D\geq4$, and $k_2=2k$}
\label{tab:drg-list}
\begin{tabular}{@{}lrrrl@{}}
\toprule
Graph & $k$ & $v$ & $D$ & Intersection array \\
\midrule
Pappus & 3 & 18 & 4 &
$\{3,2,2,1;1,1,2,3\}$ \\
Coxeter & 3 & 28 & 4 &
$\{3,2,2,1;1,1,1,2\}$ \\
Tutte $8$-cage & 3 & 30 & 4 &
$\{3,2,2,2;1,1,1,3\}$ \\
Dodecahedron & 3 & 20 & 5 &
$\{3,2,1,1,1;1,1,1,2,3\}$ \\
Desargues & 3 & 20 & 5 &
$\{3,2,2,1,1;1,1,2,2,3\}$ \\
Tutte $12$-cage & 3 & 126 & 6 &
$\{3,2,2,2,2,2;1,1,1,1,1,3\}$ \\
Biggs--Smith & 3 & 102 & 7 &
$\{3,2,2,2,1,1,1;1,1,1,1,1,1,3\}$ \\
Foster & 3 & 90 & 8 &
$\{3,2,2,2,2,1,1,1;1,1,1,1,2,2,2,3\}$ \\
\midrule
$GO(2,1)$ & 4 & 45 & 4 &
$\{4,2,2,2;1,1,1,2\}$ \\
$GD(2,1)$ & 4 & 189 & 6 &
$\{4,2,2,2,2,2;1,1,1,1,1,2\}$ \\
\midrule
$5$-cube & 5 & 32 & 5 &
$\{5,4,3,2,1;1,2,3,4,5\}$ \\
\bottomrule
\end{tabular}
\end{table}
\begin{remark*}
\begin{enumerate}
    \item Koolen and Park \cite[Theorem 20]{KP12} classified distance-regular graphs with valency $k$, diameter $D\geq 3$ and $k_2\leq \frac{3}{2}k$. 
    \item The graphs $GO(2,1)$ and $GD(2,1)$ are the line graphs of the Tutte 8-cage and Tutte 12-cage, respectively. 
    And for $k=4$, the $4$-cube is a Hadamard graph, and hence is not included in Table \ref{tab:drg-list}. 
\end{enumerate}
\end{remark*}

Then we classify $2$-walk-regular graphs with valency $k$, diameter $D\geq 4$ and $v\leq 4k+4$ vertices. 
\begin{theorem}\label{thm:4k4}
Let $\G$ be a $2$-walk-regular graph with valency $k\geq3$, diameter $D\geq4$, and $v$ vertices. Then $v\leq4k+4$ if and only if $\G$ is one of the following graphs:
\begin{enumerate}
\item a Hadamard graph with $k$ even and $v=4k$;
\item the incidence graph of a $\GDD(2,k+1;k;0,\frac{k-1}{2})$ with $k$ odd and $v=4k+4$. 
\end{enumerate}
\end{theorem}

\begin{remark*}
    Qiao, Koolen and Park \cite[Theorem 1.2]{QKP16} gave a characterization of $2$-walk-regular graphs with valency $k$, diameter $D\geq 3$ and $v\leq \alpha k$ vertices for a fixed $\alpha>2$. 
    They also showed that a $2$-walk-regular graph $\G$ with valency $k$, diameter $D\geq 3$ and $v$ vertices is distance-regular if $v\leq 3k+1$, or $v\leq 4k+1$ and $\G$ is bipartite (see \cite[Corollary 4.2]{QKP16}). They also gave an infinite family of the incidence graphs of  $\GDD(2,k+1;k;0,\frac{k-1}{2})$. Such graphs were also studied by Jin, Koolen and Lv \cite{JKL26}.
\end{remark*}

This paper is organized as follows: In Section \ref{sec:2}, we give some definitions and preliminary results. 
In Section \ref{sec:3}, we give a proof of Theorem \ref{thm:bip}. 
In Section \ref{sec:4}, we determine feasible parameters $(k,a_1)$ for $2$-walk-regular graphs with $k_2\leq 2k$ and $c_2=1$ (Lemma \ref{lem:c2-1}), and characterize the corresponding $2$-walk-regular graphs with valency four (Lemma \ref{lem:line}). 
Then we give proofs of Theorems \ref{thm:2k} and \ref{thm:drg}. 
In Section \ref{sec:5}, we determine $2$-walk-regular graphs with $v\leq 4k+4$ and $c_2=1$ (Lemma \ref{lem:v-c2}). Then we give a proof of Theorem \ref{thm:4k4}. 

\section{Definitions and preliminaries}\label{sec:2}
All graphs considered in this paper are finite, undirected and simple. 
The reader is referred to \cite{BCN89,DKT16} for more information. 

Let $\G=(V(\G),E(\G))$ be a connected graph. 
Two vertices $x,y$ are \emph{adjacent} if $\{x,y\}$ is in $E(\G)$, denoted by $x\sim y$. 
The \emph{distance} $d(x,y)=d_\G(x,y)$ between two vertices $x,y\in V(\G)$ is the length of a shortest path between $x$ and $y$ in $\G$. 
The \emph{diameter} $D=D_\G$ is the maximum distance between any two vertices in $\G$. 
For each $x\in V(\G)$, let $\G_i(x)$ be the set of vertices at distance $i$ from $x$ ($0\leq i\leq D$). 
In addition, we define $\G_{-1}(x)=\G_{D+1}(x)=\emptyset$. 
For simplicity, we write $\G(x)=\G_1(x)$. 
For $x\in V(\G)$, the \emph{valency} of $x$ in $\G$ is $|\G(x)|$. 
A graph is \emph{regular} with valency $k$ if $|\G(x)|=k$ for all $x\in V(\G)$. 
The \emph{local graph} $\Delta_x$ of $\G$ at $x$ is the subgraph induced on $\G(x)$. 
For two vertices $x,y$ at distance two, the \emph{$\mu$-graph} $M_{xy}$ is the subgraph induced on $\G(x)\cap\G(y)$. 

Let $\G$ be a graph with diameter $D$ and let $x,y$ be vertices of $\G$ at distance $i$ $(0\leq i\leq D)$. 
Let $c_i(x,y)=|\G_{i-1}(x)\cap \G(y)|$, $a_i(x,y)=|\G_{i}(x)\cap \G(y)|$, $b_i(x,y)=|\G_{i+1}(x)\cap \G(y)|$, $k_i(x)=|\G_i(x)|$. 
The intersection numbers $c_i(x,y),a_i(x,y),b_i(x,y)$ are called well-defined if the numbers do not depend on the choice of $x$ and $y$ but only on $i$. A graph $\G$ with diameter $D$ is called \emph{$t$-partially distance-regular} if these numbers $c_i(x,y)$, $a_i(x,y)$ and $b_i(x,y)$ are well-defined for $0\leq i\leq t$ (where it is implicitly assumed that $t\leq D$). 
If this is the case, then these numbers are denoted simply by $c_i$, $a_i$ and $b_i$ $(0\leq i\leq t)$. 
The relation $k_{i-1}(x)b_{i-1}=k_i(x)c_i$ shows that $k_i=k_i(x)$ $(0\leq i\leq t)$ are also well-defined. 
If $t=D$, then the graph $\G$ is called a \emph{distance-regular graph}, and the array $\{b_0=k,b_1,\ldots,b_{D-1};c_1=1,c_2,\ldots,c_D\}$ is called the \emph{intersection array} of $\G$. 
A distance-regular graph with intersection array $\{2\mu,2\mu-1,\mu,1;1,\mu,2\mu-1,2\mu\}$ is called a \emph{Hadamard graph}.

A graph $\G$ is called \emph{$t$-walk-regular} if the number of walks of every given length $\ell$ between two vertices $x,y\in V(\G)$ depends only on the distance between $x$ and $y$, provided that $d(x,y)\leq t$ (where it is implicitly assumed that the diameter of $\G$ is at least $t$). A $t$-walk-regular graph is a $t$-partially distance-regular graph (see \cite[Proposition 3.15]{DDF11}). 

The \emph{distance-$i$ matrix} $A_i$ of $\G$ is the matrix whose rows and columns are indexed by the vertices of $\G$ and the $(x,y)$-entry is $1$ whenever $d(x,y)=i$ and $0$ otherwise. 
The \emph{adjacency matrix} of $\G$ is $A_1$, and it is simply denoted by $A_\G$ or $A$. 
The eigenvalues of $\G$ are the eigenvalues of its adjacency matrix $A$. We denote the transpose of a matrix $A$ by $A^T$. 

Let $\G$ be a $t$-walk-regular graph with $v$ vertices, and let $\theta$ be an eigenvalue of $\G$ with multiplicity $m(\theta)$. 
The matrix $E$ corresponding to the orthogonal projection onto the eigenspace associated to $\theta$ is called the \emph{minimal idempotent} associated to $\theta$. The map 
\[
\rho:V(\G)\to \mathbb R^v,\qquad x\mapsto \widehat x=\sqrt{\frac{v}{m(\theta)}}Ee_x
\]
is called the \emph{normalized representation} associated to $\theta$. 
Note that there exist constants $\alpha_i$ such that $A_i\circ E=\alpha_iA_i$ $(0\leq i\leq t)$, where $\circ$ is the entrywise product (see \cite[Theorem 3.1]{DFG09}), i.e., the $(x,y)$-entry of $E$ only depends on the distance between $x$ and $y$ if $d(x,y)=i\leq t$. 
Then $\alpha_0=m(\theta)/v$, $\Vert \widehat x\Vert^2=1$, and $u_i(\theta):=\frac{\alpha_i}{\alpha_0}=\langle \widehat x,\widehat y\rangle$ $(0\leq i\leq t)$ is well-defined, where $x,y$ are any two vertices of $\G$ at distance $i$. As $AE=\theta E$, we have 
$\theta\widehat x=\sum_{y\in\G(x)}\widehat y$ and
\begin{equation}\label{eq:std}
 c_i u_{i-1}(\theta)+a_i u_i(\theta)+b_i u_{i+1}(\theta)=\theta u_i(\theta)\qquad (1\leq i\leq t-1).
\end{equation}
The sequence $(u_0(\theta),u_1(\theta),\ldots,u_t(\theta))$ is called the \emph{$t$-partial standard sequence} of $\G$ corresponding to the eigenvalue $\theta$. For a distance-regular graph with diameter $D$, the standard sequence is the $D$-partial standard sequence (see \cite[p. 128]{BCN89}). Note that for a $2$-walk-regular graph, by \eqref{eq:std}, we have
\begin{equation}\label{eq:cos}
 u_0(\theta)=1,\qquad u_1(\theta)=\frac\theta k,\qquad
 u_2(\theta)=\frac{\theta^2-a_1\theta-k}{kb_1}.
\end{equation}

A partition $ \Pi=\{B_1,B_2,\ldots,B_n\}$ of $V(\Gamma)$
is \emph{equitable} if there exist non-negative integers
$q_{ij}$ ($1\leq i,j\leq n$) such that every vertex in
$B_i$ has exactly $q_{ij}$ neighbors in $B_j$.
When $\Pi$ is equitable, the \emph{distribution diagram}
of $\Gamma$ with respect to $\Pi$ is a diagram
in which each block $B_i$ is represented by a balloon.
We draw a line from the balloon representing $B_i$
to that of $B_j$ if $q_{ij}>0$, placing the number
$q_{ij}$ above the line near $B_i$.
Inside the balloon we write $|B_i|$, and below
it we write $q_{ii}$.
If $q_{ii}=0$, we write $-$ instead of $0$.

An \emph{incidence structure} $\mathcal I=(\mathcal P,\mathcal B,I)$ consists of a set
$\mathcal P$ of points, a set $\mathcal B$ of blocks (disjoint from $\mathcal P$),
and a relation $I\subseteq\mathcal P\times\mathcal B$ called \emph{incidence}.
If $(p,B)\in I$, then we say that the point $p$ and the block $B$ are \emph{incident}.
We regard each block $B$ as a subset of $\mathcal P$.
Given an incidence structure $\mathcal I=(\mathcal P,\mathcal B,I)$, \emph{its dual
incidence structure} is defined as $\mathcal I^*=(\mathcal B,\mathcal P,I^*)$, where
$I^*=\{(B,p)\mid (p,B)\in I\}$.
The \emph{point-block incidence graph} $\Gamma(\mathcal I)$ of $\mathcal I$ is the graph
with vertex set $\mathcal P\cup\mathcal B$, where two vertices are adjacent if and
only if they are incident.
Note that the point-block incidence graph of an incidence structure is bipartite.

A \emph{group divisible design} $\mathcal D$ with parameters
$(n,m;k;\lambda_1,\lambda_2)$, denoted by
$\operatorname{GDD}(n,m;k;\lambda_1,\lambda_2)$,
is an ordered triple $(\mathcal P,\mathcal G,\mathcal B)$, where $\mathcal P$ is a set of points,
$\mathcal G$ is a partition of $\mathcal P$ into $m$ subsets of size $n$
(called groups), and $\mathcal B$ is a collection of $k$-subsets of $\mathcal P$
(called blocks), such that each pair of points from the same group appears in exactly
$\lambda_1$ blocks, and each pair of points from different groups appears in exactly
$\lambda_2$ blocks. 
The triple $\mathcal I=(\mathcal P,\mathcal B,I)$, with the natural incidence relation $I$,
forms an incidence structure, and we consider its dual
$\mathcal I^*=(\mathcal B,\mathcal P,I^*)$.
If there exists a partition $\mathcal G'$ of $\mathcal B$ such that the triple
$(\mathcal B,\mathcal G',\mathcal P)$ is a
$\operatorname{GDD}(n,m;k;\lambda_1,\lambda_2)$, then we say that $\mathcal D$
is a \emph{group divisible design with the dual property} with parameters
$(n,m;k;\lambda_1,\lambda_2)$, denoted by
$\operatorname{GDDDP}(n,m;k;\lambda_1,\lambda_2)$. 

For the convenience of the reader, we list the following known results. 
\begin{lemma}[{cf. \cite[Proposition 4.1.6]{BCN89}}]
    \label{lem:bc}
    Let $\G$ be a $t$-walk-regular graph with diameter $D\geq t$. 
    If $i,j\leq t$ and $i+j\leq D$, then $b_i\geq c_j$ holds. 
\end{lemma}

\begin{lemma}[{cf. {\cite[Proposition~5.1]{CDK13}}}]
    \label{lem:local-eig}
    Let $\G$ be a $2$-walk-regular graph with diameter $D\geq 3$ and distinct eigenvalues $k=\theta_0>\theta_1>\cdots>\theta_d$. 
    Let $\Delta$ be a local graph of $\G$ with smallest eigenvalue $\eta$. Then 
    \[
    \eta\geq -1-\frac{b_1}{\theta_1+1}. 
    \]
\end{lemma}

\begin{lemma}[{cf. \cite[Lemma~2.5]{QPK19}}]
\label{lem:quad-eig}
    Let $\G$ be a $2$-walk-regular graph with diameter $D\geq 3$ and distinct eigenvalues $k=\theta_0>\theta_1>\cdots>\theta_d$. 
    If $\G$ contains an induced quadrangle, then
    \[
    \theta_1\leq b_1-1. 
    \]
\end{lemma}

\begin{lemma}[{cf. \cite[Theorem~1.3]{QPK19}}]
    \label{lem:large-c2}
    Let $\G$ be a $2$-walk-regular graph with valency $k\geq3$ and diameter $D\geq4$. 
    If $c_2>k/3$, then $\G$ is a Hadamard graph,
    the $5$-cube, or the incidence graph of a
    $\GDD(2,m;k;0,c_2)$ for some integer $m$.
\end{lemma}

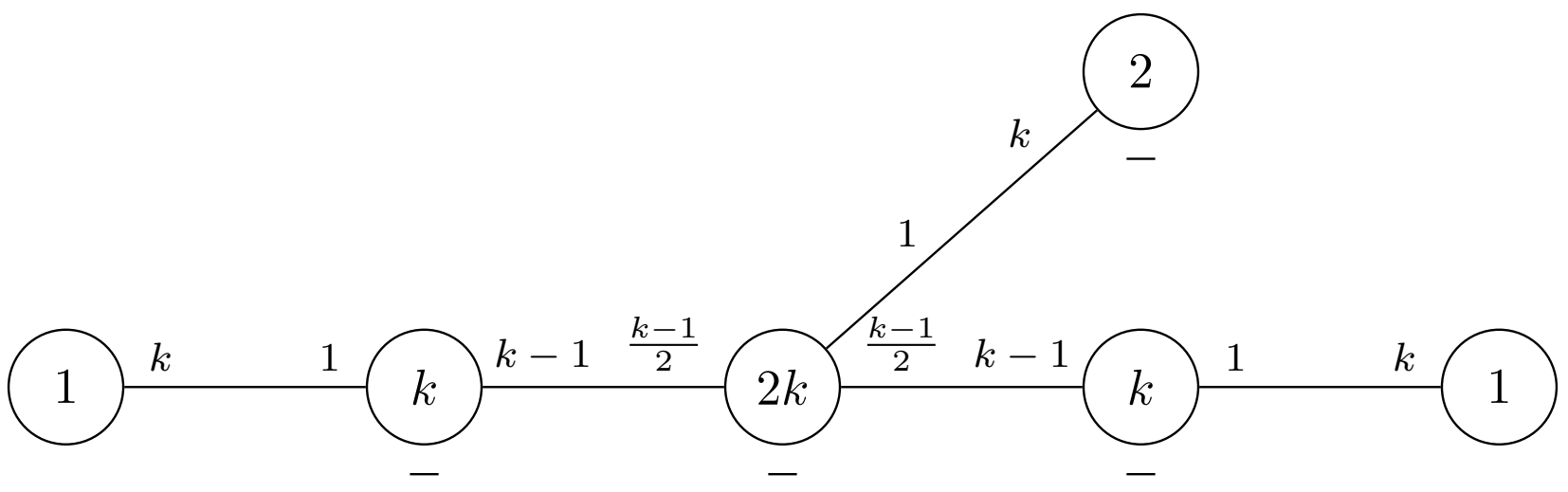
\begin{figure}[htbp]
\centering
\begin{tikzpicture}[
  balloon/.style={circle,draw,line width=.45pt,
                  minimum size=8mm,inner sep=1pt},
  every node/.style={font=\small},
  edge label/.style={font=\scriptsize,inner sep=1pt},line width=.45pt]
\node[balloon] (x)  at (-5,0)   {$1$};
\node[balloon] (n)  at (-2.5,0) {$k$};
\node[balloon] (s)  at (0,0)    {$2k$};
\node[balloon] (np) at (2.5,0)  {$k$};
\node[balloon] (xp) at (5,0)    {$1$};
\node[balloon] (u)  at (2.5,2.2) {$2$};
\draw (x) -- node[edge label,pos=.15,above=2pt] {$k$}
             node[edge label,pos=.85,above=2pt] {$1$} (n);
\draw (n) -- node[edge label,pos=.25,above=2pt] {$k-1$}
             node[edge label,pos=.75,above=2pt] {$\tfrac{k-1}{2}$} (s);
\draw (s) -- node[edge label,pos=.25,above=2pt] {$\tfrac{k-1}{2}$}
             node[edge label,pos=.75,above=2pt] {$k-1$} (np);
\draw (np) -- node[edge label,pos=.15,above=2pt] {$1$}
              node[edge label,pos=.85,above=2pt] {$k$} (xp);
\draw (s) -- node[edge label,pos=.38,above left=1pt] {$1$}
             node[edge label,pos=.80,above left=1pt] {$k$} (u);
\foreach \p in {n,s,np,u} {
  \draw ([xshift=-1mm,yshift=-2mm]\p.south) --
        ([xshift=1mm,yshift=-2mm]\p.south);
}
\end{tikzpicture}
\caption{Distribution diagram of the incidence graph of $\GDD(2,k+1;k;0,\frac{k-1}{2})$}
\label{fig:gddk}
\end{figure}

\begin{lemma}[{cf. \cite[Theorem~2]{QDK16}}]
\label{lem:gdd}
Let $k\geq3$ be odd, and let $\mathcal D$ be a $\GDD(2,k+1;k;0,\frac{k-1}{2})$ with incidence graph $\G$. 
Then $\G$ is $2$-walk-regular with $v=4k+4$ vertices. 
Moreover, for each vertex $x$, there exists a unique vertex $x'$ at distance four from $x$, 
         and the partition 
         \[
         \Pi(x)=\{\{x\},\G(x),\G_2(x),\G(x'), \G_3(x)\setminus\G(x'),\{x'\}\}
         \]
        is equitable with distribution diagram as in Figure~\ref{fig:gddk}.
\end{lemma}

\begin{lemma}[{cf. \cite[Propositions 3.6 and 3.4]{CDK13}}]
\label{lem:walk-transfer}
Let $\G$ be a graph with valency $k\geq 3$, diameter $D\geq 3$ and girth $g$. 
\begin{enumerate}
    \item If $\G$ is $1$-walk-regular and $k=3$, then $\G$ is $2$-walk-regular. 
    \item If $\G$ is $3$-walk-regular and $g\geq 6$, then $L(\G)$ is $2$-walk-regular. 
\end{enumerate}
\end{lemma}

\begin{lemma}[{cf. \cite[Lemma~6]{KP10}}]
\label{lem:theta1}
Let $\G$ be a $2$-walk-regular graph with valency $k\geq3$,
diameter $D\geq4$, and distinct eigenvalues $ k=\theta_0>\theta_1>\cdots>\theta_d$. 
Then
\[
 \theta_1\geq\frac{a_1+\sqrt{a_1^2+4k}}{2}.
\]
\end{lemma}
\begin{proof}
    Let $x,y\in V(\G)$ with $d(x,y)=4$. Then the graph induced on $\{x\}\cup \G(x)\cup\{y\}\cup \G(y)$ has largest eigenvalue $\frac{a_1+\sqrt{a_1^2+4k}}{2}$ with multiplicity 2. The result follows by interlacing. 
\end{proof}

\section{Proof of Theorem \ref{thm:bip}}\label{sec:3}

\begin{proof}[Proof of Theorem \ref{thm:bip}]
	Let $\G$ be a $2$-walk-regular graph with valency $k\geq3$, diameter $D\geq4$,
	satisfying $k_2\leq2k$ and $c_2\geq2$.

	Let $\tau>0$ and $-\sigma<0$ be the roots of $t^2-a_1t-k$. Then
	\begin{equation}\label{eq:roots}
		k=\sigma(a_1+\sigma),\qquad
		\tau=a_1+\sigma,\qquad
		b_1=(a_1+\sigma+1)(\sigma-1),\qquad
		\sigma>1.
	\end{equation}

	Let $\Delta$ be a local graph, and let $\mu_\Delta$ be the average number
	of common neighbors in $\Delta$ over its unordered pairs of nonadjacent vertices,
	that is 
	\[
		\mu_\Delta
		=\frac{2}{kb_1}
		\sum_{\substack{\{y,z\}\subseteq V(\Delta)\\ y\not\sim z}}
		|\Delta(y)\cap\Delta(z)|.
	\]
	For any two vertices $p,q$ of $\G$ at distance 2, let $\delta_{pq}$ be the average degree
	of the $\mu$-graph $M=M_{pq}$, that is 
	\[
		\delta_{pq}
		=\frac{1}{c_2}\sum_{y\in V(M)}|M(y)|.
	\]

	\noindent\textbf{Claim 1.} We have $\theta_1\geq\tau$, and the local graph $\Delta$
	has smallest eigenvalue at least $-\sigma$. Moreover,
	\begin{equation}\label{eq:mu-upper}
		\mu_\Delta\leq\frac{a_1}{\sigma+1},\qquad
		\mu_\Delta\leq\frac{a_1(a_1-1)}{b_1}.
	\end{equation}

	\begin{proof}[Proof of Claim 1.]
		By Lemma \ref{lem:theta1}, we have $\theta_1\geq\tau$.
		Thus Lemma \ref{lem:local-eig} shows that the smallest eigenvalue
		of the local graph $\Delta$ is at least
		\[
			-1-\frac{b_1}{\theta_1+1}
			\geq-1-\frac{b_1}{\tau+1}
			=-\sigma.
		\]

		In the local graph $\Delta$, there are $kb_1$ ordered nonadjacent pairs
		of vertices and $ka_1(a_1-1)$ walks of length two with distinct endpoints,
		among which $\tr(A_\Delta^3)$ have adjacent endpoints. Hence
		\begin{equation}\label{eq:trace}
			kb_1\mu_\Delta
			=ka_1(a_1-1)-\tr(A_\Delta^3).
		\end{equation}
		Since $\tr(A_\Delta^3)\geq0$, this shows the second inequality
		in \eqref{eq:mu-upper}.

		Let $\rho=a_1/(\sigma+1)$.
		Since the spectrum of $A_\Delta$ lies in $[-\sigma,a_1]$,
		the polynomial $f(t)=(t+\sigma)(t-\rho)^2$ is nonnegative on the spectrum.
		Since $a_1$ is an eigenvalue of $A_\Delta$, we have
		$\tr f(A_\Delta)\geq f(a_1)$.
		Using $\tr(A_\Delta)=0$ and $\tr(A_\Delta^2)=ka_1$, this gives
		\[
			\tr(A_\Delta^3)+(\sigma-2\rho)ka_1+\sigma\rho^2k
			\geq(a_1+\sigma)(a_1-\rho)^2
			=\sigma\rho^2k.
		\]
		Hence $\tr(A_\Delta^3)\geq ka_1(2\rho-\sigma)$.
		Substituting this bound into \eqref{eq:trace}, we obtain
		\[
			\mu_\Delta
			\leq\frac{a_1(a_1-1-2\rho+\sigma)}{b_1}
			=\frac{a_1}{\sigma+1}.
		\]
		This shows the first inequality in \eqref{eq:mu-upper}.
	\end{proof}

	\noindent\textbf{Claim 2.} $\G$ contains an induced quadrangle,
	$\tau\leq\theta_1\leq b_1-1$ and $\sigma\geq2$.

	\begin{proof}[Proof of Claim 2.]
		Assume that $\G$ contains no induced quadrangle.
		Then every $\mu$-graph of $\G$ is a clique.
		Hence in the local graph $\Delta$, two nonadjacent vertices have
		exactly $c_2-1$ common neighbors.
		Claim 1 gives
		\begin{equation}\label{eq:c2-mu}
			c_2-1=\mu_\Delta\leq\frac{a_1}{\sigma+1}.
		\end{equation}
		Together with $b_1\leq2c_2$ and \eqref{eq:roots}, we have
		\[
			(a_1+\sigma+1)(\sigma^2-1)
			=b_1(\sigma+1)
			\leq2c_2(\sigma+1)
			\leq2(a_1+\sigma+1),
		\]
		and hence $\sigma\leq\sqrt3<2$.
		Moreover, \eqref{eq:c2-mu} implies
		$a_1\geq(c_2-1)(\sigma+1)>2$, since $c_2\geq2$ and $\sigma>1$.

		The local graph $\Delta$ is connected, since its nonadjacent vertices
		have $c_2-1>0$ common neighbors.
		By Claim 1, the smallest eigenvalue of $\Delta$ is at least $-\sigma>-2$.
		By \cite[Corollary 3.12.3]{BCN89}, it must be an
		odd cycle, contrary to $a_1>2$.
		Hence $\G$ contains an induced quadrangle.

		Now Claim 1 and Lemma \ref{lem:quad-eig} give
		$\tau\leq\theta_1\leq b_1-1$.
		Since $b_1=(\tau+1)(\sigma-1)$, we have $\sigma\geq2$.
	\end{proof}

	\noindent\textbf{Claim 3.}
	For any two vertices $p,q\in V(\G)$ with $d(p,q)=2$,
	\begin{equation}\label{eq:delta-lb}
		\delta_{pq}\geq\frac{2c_2}{\sigma}-\sigma.
	\end{equation}
	The average $\overline{\delta}$ of $\delta_{pq}$ over all unordered
	pairs of vertices of $\G$ at distance 2 equals the average of $\mu_\Delta$
	over all local graphs of $\G$.

	\begin{proof}[Proof of Claim 3.]
		By Claim 2, we have
		$\tau\leq\theta_1\leq b_1-1$.
		Write $u_1=u_1(\theta_1)$ and $u_2=u_2(\theta_1)$.
		Equation \eqref{eq:cos} gives
		\[
			u_2
			=\frac{(\theta_1-\tau)(\theta_1+\sigma)}{kb_1}
			\geq0,\qquad
			u_1-u_2
			=\frac{(k-\theta_1)(\theta_1+1)}{kb_1}
			>0,
		\]
		\[
			\eta
			:=\frac{1-u_2}{u_1-u_2}
			=1+\frac{b_1}{\theta_1+1}
			\leq1+\frac{b_1}{\tau+1}
			=\sigma,\qquad
			2\leq\eta\leq\sigma.
		\]

		Consider a normalized representation $z\mapsto\widehat z$
		associated to $\theta_1$.
		Fix $p,q\in V(\G)$ with $d(p,q)=2$.
		Let $M=M_{pq}$ and
		\[
			Z=\sum_{z\in V(M)}\widehat z.
		\]
		Since vertices of $M$ are pairwise at distance one or two in $\G$,
		\[
			\begin{aligned}
				\|Z\|^2
				&=c_2(1-u_2)
				+c_2\delta_{pq}(u_1-u_2)
				+c_2^2u_2,\\
				\langle Z,\widehat p+\widehat q\rangle
				&=2c_2u_1,
				&\|\widehat p+\widehat q\|^2
				&=2(1+u_2).
			\end{aligned}
		\]
		Applying Cauchy--Schwarz to $Z$ and $\widehat p+\widehat q$ gives
		\[
			(2c_2u_1)^2
			\leq2(1+u_2)
			\left(
				c_2(1-u_2)
				+c_2\delta_{pq}(u_1-u_2)
				+c_2^2u_2
			\right).
		\]
		Since $1+u_2>0$ and $u_1-u_2>0$, this is equivalent to
		\[
			\delta_{pq}
			\geq
			\frac{
				c_2(2u_1^2/(1+u_2)-u_2)-(1-u_2)
			}{
				u_1-u_2
			}
			=c_2-\eta-\frac{c_2u_1(\eta-2)}{1+u_2}.
		\]
		Moreover, using \eqref{eq:roots},
		\[
			0\leq u_1(\eta-2)
			=\frac{\sigma-2}{\sigma}
			-\frac{
				(\theta_1-\tau)(\theta_1+2-\sigma)
			}{
				k(\theta_1+1)
			}
			\leq\frac{\sigma-2}{\sigma},
		\]
		since $\theta_1\geq\tau$ and
		$\theta_1+2-\sigma\geq a_1+2>0$.
		As $u_2\geq0$ and $\eta\leq\sigma$, it follows that
		\[
			\delta_{pq}
			\geq
			c_2-\sigma-\frac{c_2(\sigma-2)}{\sigma}
			=\frac{2c_2}{\sigma}-\sigma.
		\]

		Let $\mathcal S_2$ be the set of unordered distance-two pairs of $V(\G)$.
		We count the number of induced diamonds (see Figure \ref{fig:diamond}) in $\G$.
		Then
		\[
			\sum_{\{p,q\}\in\mathcal S_2}
			\left(\frac{1}{2}c_2\delta_{pq}\right)
			=
			\frac{1}{2}
			\sum_{x\in V(\G)}
			\left(\frac{1}{2}kb_1\cdot\mu_{\Delta_x}\right),
		\]
		where the $\mu$-graph $M_{pq}$ has
		$\frac{1}{2}c_2\delta_{pq}$ edges, and the local graph $\Delta_x$
		of $\G$ at vertex $x$ has $\frac{1}{2}kb_1$ unordered pairs
		of nonadjacent vertices.
		Thus we have
		\[
			\frac{1}{|\mathcal S_2|}
			\sum_{\{p,q\}\in\mathcal S_2}\delta_{pq}
			=
			\frac{1}{v}
			\sum_{x\in V(\G)}\mu_{\Delta_x},
		\]
		where
		\[
			|\mathcal S_2|
			=\frac{1}{2}vk_2
			=\frac{vkb_1}{2c_2}.
		\]

		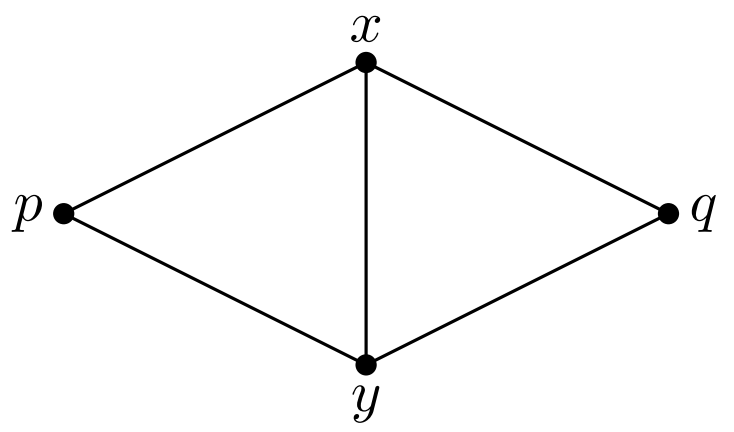
\begin{figure}[htbp]
			\centering
			\begin{tikzpicture}[line width=.6pt]
				\coordinate (p) at (-2,0);
				\coordinate (q) at (2,0);
				\coordinate (x) at (0,1);
				\coordinate (y) at (0,-1);

				\fill (p) circle (2pt) node[left] {$p$};
				\fill (q) circle (2pt) node[right] {$q$};
				\fill (x) circle (2pt) node[above] {$x$};
				\fill (y) circle (2pt) node[below] {$y$};

				\draw (p) -- (x) -- (q) -- (y) -- (p);
				\draw (x) -- (y);
			\end{tikzpicture}
			\caption{Induced diamond}
			\label{fig:diamond}
		\end{figure}
	\end{proof}

	\noindent\textbf{Claim 4.} $a_1=0$.

	\begin{proof}[Proof of Claim 4.]
		Let $\overline\delta$ be the average of
		$\delta_{pq}$ over all unordered pairs of the vertices of $\G$ at distance 2.
		Then by Claim 1 and Claim 3, we have
		\begin{equation}\label{eq:delta-bounds}
			\frac{2c_2}{\sigma}-\sigma
			\leq\overline{\delta}
			\leq
			\min\left\{
				\frac{a_1}{\sigma+1},
				\frac{a_1(a_1-1)}{b_1}
			\right\}.
		\end{equation}
		Since $b_1\leq2c_2$ and $\sigma\geq2$ (Claim 2),
		using \eqref{eq:roots}, we have
		\[
			\frac{2c_2}{\sigma}-\sigma
			\geq\frac{b_1}{\sigma}-\sigma
			=a_1-\frac{a_1+1}{\sigma}
			\geq\frac{a_1-1}{2},\qquad
			b_1\geq a_1+3.
		\]
		By \eqref{eq:delta-bounds}, we have
		\begin{equation}\label{eq:delta-chain}
			\frac{a_1-1}{2}
			\leq\overline\delta
			\leq
			\min\left\{
				\frac{a_1}{3},
				\frac{a_1(a_1-1)}{a_1+3}
			\right\}.
		\end{equation}

        Assume $a_1\geq2$. 
		The lower bound and first upper bound in \eqref{eq:delta-chain}
		imply $a_1\leq3$.
		Then the lower bound and second upper bound in \eqref{eq:delta-chain}
		give $a_1\geq3$, hence $a_1=3$.
		By \eqref{eq:delta-chain}, we have $\overline{\delta}=1$.
		Then $\overline{\delta}\leq a_1/(\sigma+1)$ implies $\sigma=2$.
		Hence $b_1=6$ by \eqref{eq:roots}.
		By \eqref{eq:delta-bounds}, $\overline{\delta}\geq c_2-2$, so $c_2\leq3$.
		Since $b_1\leq2c_2$, we have $c_2=3$.

		Choose $p,q\in V(\G)$ with $d(p,q)=2$.
		Note that the $\mu$-graph $M_{pq}$ has 3 vertices.
		Then \eqref{eq:delta-lb} implies $\delta_{pq}\geq1$.
		It follows that $M_{pq}$ has at least two edges and
		$\delta_{pq}\geq4/3$, contradicting $\overline{\delta}=1$.
		Hence $a_1\leq1$.

		Assume $a_1=1$.
		Then every local graph is $1$-regular and every independent set
		in a local graph has at most $k/2$ vertices.
		Note that $\G$ contains no induced diamond and each $\mu$-graph
		contains no edges.
		Choose $x,y,z\in V(\G)$ with $d(x,z)=4$ and
		$d(x,y)=d(y,z)=2$.
		Then vertices in $M_{xy}$ and $M_{yz}$ form an independent set
		of $\Delta_y$ with $2c_2$ vertices.
		Hence $2c_2\leq k/2$.
		As $a_1=1$ and $b_1\leq2c_2$, we have
		$k=b_1+2\leq2c_2+2$.
		Then $c_2\leq1$, a contradiction.
		Hence $a_1=0$.
	\end{proof}

	By Claim 4, $a_1=0$ and $b_1=k-1$.
	Assume $\G$ is nonbipartite.
	Then Lemma \ref{lem:large-c2} implies $c_2\leq k/3$.
	As $b_1\leq2c_2$, we have $k-1=b_1\leq2c_2\leq 2k/3$. 
	Thus $k\leq3$ and $c_2\leq1$, contrary to $c_2\geq2$.
\end{proof}

\section{Proofs of Theorems \ref{thm:2k} and \ref{thm:drg}}\label{sec:4}

\begin{lemma}\label{lem:c2-1}
Let $\G$ be a $2$-walk-regular graph with valency $k\geq3$, diameter
$D\geq4$, and $k_2\leq2k$. If $c_2=1$, then
\[
	(k,a_1)=(3,0)\quad\text{or}\quad(4,1).
\]
In particular, if $k\geq5$, then $c_2\geq2$.
\end{lemma}

\begin{proof}
	Note that $1\leq b_1\leq2c_2=2$.
	As $c_2=1$, $\G$ contains no induced diamond and every local graph is a disjoint union of cliques of size $a_1+1$.
	Since $a_1+1=k-b_1$, we have $k-b_1\mid k$, that is $k-b_1\mid b_1$.
	It follows that $b_1=2$, and $k=3$ or $4$.
\end{proof}

\begin{lemma}\label{lem:diam}
    Let $H$ be a cubic $3$-walk-regular graph with diameter $D_H\geq 3$ and girth at least six. 
    Then $D_H\geq 4$ if and only if $D_{L(H)}\geq 4$. 
\end{lemma}
\begin{proof}
    Assume $D_H\geq 4$ and consider a path $x_0x_1x_2x_3x_4$ with $d_H(x_0,x_4)=4$. 
    Since $H$ is $3$-walk-regular, we have $b_3\neq 0$. Take a neighbor $z$ of $x_4$ with $d_H(x_1,z)=4$. Then $d_{L(H)}(x_0x_1,x_4z)=4$. 
    It follows that $D_{L(H)}\geq 4$. 

    If $D_H=3$, then $H$ is distance-regular, which is the Heawood graph (see \cite[Theorem 7.5.1]{BCN89}). 
    Since $H$ is bipartite with diameter three, any two edges have
	endpoints in the same bipartition class at distance at most two.
	Hence $D_{L(H)}\leq 3$.
\end{proof}

\begin{lemma}\label{lem:line}
\begin{enumerate}
	\item For a $2$-walk-regular graph $\G$ with valency $k=4$, diameter $D\geq4$ and $(a_1,c_2)=(1,1)$, there exists a unique cubic $3$-walk-regular graph $H$ with diameter $D_H\geq 4$ and girth $g_H\geq6$, such that $\G=L(H)$.
	\item Conversely, for a cubic $3$-walk-regular graph $H$ with diameter $D_H\geq 4$ and girth $g_H\geq6$, the line graph $L(H)$ is $2$-walk-regular with valency $4$, diameter $D_{L(H)}\geq 4$ and $(a_1,c_2)=(1,1)$.
\end{enumerate}
\end{lemma}

\begin{proof}
\begin{enumerate}
	\item Note that the local graph of $\G$ consists of two disjoint edges.
	Let $V(H)$ be the set of triangles of $\G$, and two vertices of $H$ are adjacent if the triangles share a common vertex.
	Then $H$ is cubic.
	Let $N$ be the triangle-vertex incidence matrix of $\G$, which is also the vertex-edge incidence matrix of $H$.
	Then $A_\G=N^TN-2I$ implies $\G=L(H)$. 
    
	A triangle in $H$ would give an extra triangle in $\G$, and an induced quadrangle in $H$ would give an induced quadrangle in $\G$.
	Hence $g_H\geq5$.
	Since $\G=L(H)$ and $D_H+1\geq D_{L(H)}\geq4$, we have $D_H\geq3$. 
    
	Choose a path $x_0x_1x_2x_3$ in $H$ with $d_H(x_0,x_3)=3$.
	Let $x_2x_4$ and $x_3x_5$ be extra edges.
	Note that $d_\G(x_0x_1,x_2x_3)=2$ and $d_\G(x_0x_1,x_2x_4)=2$. 
    If $d_\G(x_0x_1,x_3x_5)=2$, then $x_1\sim x_5$, contradicting $g_H\geq 5$.
	Hence $a_2=1$ in $\G$.
	If $H$ contains an induced pentagon $w_0w_1w_2w_3w_4w_0$, then $d_\G(w_0w_1,w_2w_3)=2$ and $w_3w_4,w_2w_6\in \G_2(w_0w_1)$, where $w_2w_6$ is an extra edge, contrary to $a_2=1$. 
	This shows $g_H\geq6$.

	Consider the vertex-edge incidence matrix of $H$.
	Let
	\[
		P=NN^T=A_H+3I,\qquad Q=N^TN=A_\G+2I.
	\]
	For every $\ell\geq1$, we have
	\[
		P^\ell=NQ^{\ell-1}N^T,\qquad
		Q^{\ell+1}=N^TP^\ell N.
	\]
	By the first identity, the diagonal of $P^\ell$ is a constant $\beta_0$, as the three edges incident with a vertex of $H$ form a triangle in $\G$.
	Let the entries of $Q^{\ell+1}$ at distance $i\leq2$ in $\G$ be $\gamma_i$.
	Then we successively determine entries $\beta_i$ of $P^\ell$ at distance $i=1,2,3$ in $H$ through
	\[
		2\beta_0+2\beta_1=\gamma_0,\qquad
		\beta_0+2\beta_1+\beta_2=\gamma_1,\qquad
		\beta_1+2\beta_2+\beta_3=\gamma_2.
	\]
	Since powers of $A_H=P-3I$ are linear combinations of powers of $P$, the graph $H$ is $3$-walk-regular. And Lemma \ref{lem:diam} implies $D_H\geq 4$.  

	Assume $\G=L(H)$ with $H$ cubic and girth $g_H\geq6$.
	Then each vertex of $H$ is incident with three edges.
	Hence each vertex of $H$ corresponds to a triangle in $\G$.
	As $H$ is triangle-free, every triangle of $\G$ arises in this way.
	And two vertices of $H$ are adjacent if the triangles in $\G$ share a common vertex.
	This shows the uniqueness of $H$.

	\item Conversely, Lemma \ref{lem:walk-transfer} implies that $L(H)$ is $2$-walk-regular, and Lemma \ref{lem:diam} implies $D_{L(H)}\geq 4$. 
	Since $g_H\geq6$, the local graph of $L(H)$ consists of two disjoint edges and $a_1=1$.
	Then $L(H)$ contains no induced quadrangles and $c_2=1$.
\end{enumerate}
\end{proof}

\begin{proof}[Proof of Theorem~\ref{thm:2k}]
Assume $c_2\geq 2$. Note that $b_1\leq 2c_2$. 
Then $\G$ is bipartite by Theorem~\ref{thm:bip}.
We have $b_1=k-1$ and $b_2=k-c_2$. 
Since $b_1\leq 2c_2$ and $c_2\leq b_2$ (Lemma~\ref{lem:bc}), we have \[
(k-1)/2\leq c_2\leq k/2.
\]

If $k$ is even, then $c_2=k/2$ and $k\geq 4$. 
Then $\G$ is a Hadamard graph by \cite[Theorem~1.9.3]{BCN89}, giving \eqref{item:had}. 

If $k$ is odd, then $c_2=(k-1)/2$ and $k\geq 5$. 
Hence $c_2>k/3$. 
Lemma~\ref{lem:large-c2} shows that 
$\G$ is the $5$-cube, giving \eqref{item:q5}, or the
incidence graph of a $\GDD(2,m;k;0,c_2)$. 
In the latter case, the vertices at distance two from a point
are precisely the $2(m-1)$ points outside its group. Hence $2(m-1)=k_2=\frac{kb_1}{c_2}=2k$, that is $m=k+1$. 
Lemma~\ref{lem:gdd} shows  that $D=4$, giving \eqref{item:gdd}. 

If $c_2=1$, Lemma~\ref{lem:c2-1} gives
$(k,a_1)=(3,0)$ or $(4,1)$.
In the cubic case, $a_1=0$ and $c_2=1$ show that the girth is at least five, giving \eqref{item:cubic}. 
In the quartic case, Lemma~\ref{lem:line} gives \eqref{item:line}.

Conversely, the $5$-cube and Hadamard graphs are distance-regular, hence $2$-walk-regular, with $k_2=2k$ and $k_2=2k-2$, respectively. 
The graphs in \eqref{item:gdd} are $2$-walk-regular 
with $k_2=2k$ by Lemma~\ref{lem:gdd}. 
For the graphs in \eqref{item:cubic}, Lemma~\ref{lem:walk-transfer} gives
$2$-walk-regularity. 
The girth condition implies
$(a_1,c_2)=(0,1)$, so $k_2=6=2k$.
For the graphs in \eqref{item:line}, Lemma~\ref{lem:line}
gives $2$-walk-regularity and $(k,a_1,c_2)=(4,1,1)$,
so $k_2=8=2k$.
\end{proof}

\begin{proof}[Proof of Theorem~\ref{thm:drg}]
Apply Theorem~\ref{thm:2k} to distance-regular
graphs with $k_2\leq2k$. 
The only possible cases are \eqref{item:had}, \eqref{item:q5}, \eqref{item:cubic} and \eqref{item:line}. 
The result follows from the classifications of distance-regular graphs of valency $3$ or $4$ (see \cite[Theorem 7.5.1]{BCN89} and \cite[Theorem 1.1]{BK99}). 
\end{proof}

\section{Proof of Theorem \ref{thm:4k4}}\label{sec:5}

\begin{lemma}\label{lem:v-2k}
Let $\G$ be a $2$-walk-regular graph with valency
$k\geq3$, diameter $D\geq4$, and $v\leq 4k+4$ vertices.
Then $k_2\leq 2k$.
\end{lemma}
\begin{proof}
Choose two vertices $x,y$ with $d(x,y)=4$. Note that  $\{x\},\G(x),\G_2(x),\G(y),\{y\}$ are pairwise disjoint. 
Hence we have 
\begin{equation}\label{eq:v}
    2+2k+k_2\leq v\leq 4k+4,
\end{equation}
that is $k_2 \leq 2k+2$. 
Suppose $k_2>2k$. Then $2c_2+1\leq b_1\leq k-1$, and
\[
 k_2=\frac{kb_1}{c_2}\geq2k+\frac{k}{c_2}>2k+2,
\]
a contradiction.
\end{proof}

\begin{lemma}\label{lem:v-c2}
Let $\G$ be a $2$-walk-regular graph with valency
$k\geq3$, diameter $D\geq4$, and $v\leq 4k+4$ vertices.
If $c_2=1$, then $k=3$ and $\G$ is the incidence graph of a
$\GDD(2,4;3;0,1)$.
\end{lemma}
\begin{proof}
Lemma \ref{lem:v-2k} gives $k_2\leq 2k$. Then we only need to consider \eqref{item:cubic} and \eqref{item:line} of Theorem \ref{thm:2k}. 

In \eqref{item:line} of Theorem \ref{thm:2k}, $k=4$ and $\G=L(H)$, where $H$ is a
cubic $3$-walk-regular graph with diameter at least four and girth at least six. 
We have $|H_2(x)|=6$, and $|V(H)|\geq 14$ by \eqref{eq:v}. 
Then 
\[
v=|E(H)|=\frac{3}{2}|V(H)|\geq 21>4k+4.
\]

In \eqref{item:cubic} of Theorem \ref{thm:2k}, we have $k=3$ and girth at least five. 
Then $a_1=0$, $c_2=1$ and $k_2=6$. 
Note that $kv=2|E(\G)|$, that is $2\mid v$. 
Equation \eqref{eq:v} implies $v\in\{14,16\}$. 
Lemma~\ref{lem:bc} gives
$b_2\geq c_2=1$, so $a_2=2-b_2\in\{0,1\}$.

Suppose $a_2=1$. 
For every vertex $x$, each edge in $\G_2(x)$ determines a unique pentagon through $x$. 
Thus each vertex lies in $k_2a_2/2=3a_2$ pentagons and the total number of pentagons is $\frac{3va_2}{5}$, that is $5\mid v$, which contradicts $v\in\{14,16\}$. 

It follows that $a_2=0$ and $b_2=2$. 
Take two vertices $x,y$ of $\G$ with $d(x,y)=4$. 
Note that there are at most $kb_1= 6$ edges between $\G(y)$ and $\G_2(x)$, and $k_2b_2=12$ edges between $\G_2(x)$ and $\G_3(x)$. 
It follows that $|\G_3(x)\setminus \G(y)|\geq \frac{12-6}{3}=2$ and $v\geq 1+k+k_2+(k+2)+1=16$. 
Then we have $v=16$, $\G(y)\subseteq\G_3(x)$ and $|\G_3(x)\setminus\G(y)|=2$. 
Note that there are six edges between $\G_2(x)$ and each of $\G(y)$ and $\G_3(x)\setminus \G(y)$. 
We have $c_3(x,w)=2$ for $w\in \G(y)$, and $c_3(x,z)=3$ for $z\in \G_3(x)\setminus \G(y)$. 
Since $c_2=1$, each vertex in $\G_2(x)$ has exactly one neighbor in each of $\G(y)$ and $\G_3(x)\setminus\G(y)$.  
Then the partition $\Pi(x)=\{\{x\},\G(x),\G_2(x),\G(y), \G_3(x)\setminus\G(y),\{y\}\}$ has the distribution diagram as in Figure \ref{fig:gdd3}. 
Since $\G$ is bipartite with $v=16$ and $k_2=6$, every vertex has a unique vertex at distance four. 
Hence $\G$ is the incidence graph of $\GDD(2,4;3;0,1)$. 

\begin{figure}[htbp]
\centering
\begin{tikzpicture}[
  balloon/.style={circle,draw,line width=.45pt,
                  minimum size=8mm,inner sep=1pt},
  every node/.style={font=\small},
  edge label/.style={font=\scriptsize,inner sep=1pt},line width=.45pt]
\node[balloon] (x)  at (-5,0)   {$1$};
\node[balloon] (n)  at (-2.5,0) {$3$};
\node[balloon] (s)  at (0,0)    {$6$};
\node[balloon] (np) at (2.5,0)  {$3$};
\node[balloon] (xp) at (5,0)    {$1$};
\node[balloon] (u)  at (2.5,2.2) {$2$};
\draw (x) -- node[edge label,pos=.15,above=2pt] {$3$}
             node[edge label,pos=.85,above=2pt] {$1$} (n);
\draw (n) -- node[edge label,pos=.25,above=2pt] {$2$}
             node[edge label,pos=.75,above=2pt] {$1$} (s);
\draw (s) -- node[edge label,pos=.25,above=2pt] {$1$}
             node[edge label,pos=.75,above=2pt] {$2$} (np);
\draw (np) -- node[edge label,pos=.15,above=2pt] {$1$}
              node[edge label,pos=.85,above=2pt] {$3$} (xp);
\draw (s) -- node[edge label,pos=.38,above left=1pt] {$1$}
             node[edge label,pos=.80,above left=1pt] {$3$} (u);
\foreach \p in {n,s,np,u} {
  \draw ([xshift=-1mm,yshift=-2mm]\p.south) --
        ([xshift=1mm,yshift=-2mm]\p.south);
}
\end{tikzpicture}
\caption{Distribution diagram of the incidence graph of $\GDD(2,4;3;0,1)$}
\label{fig:gdd3}
\end{figure}
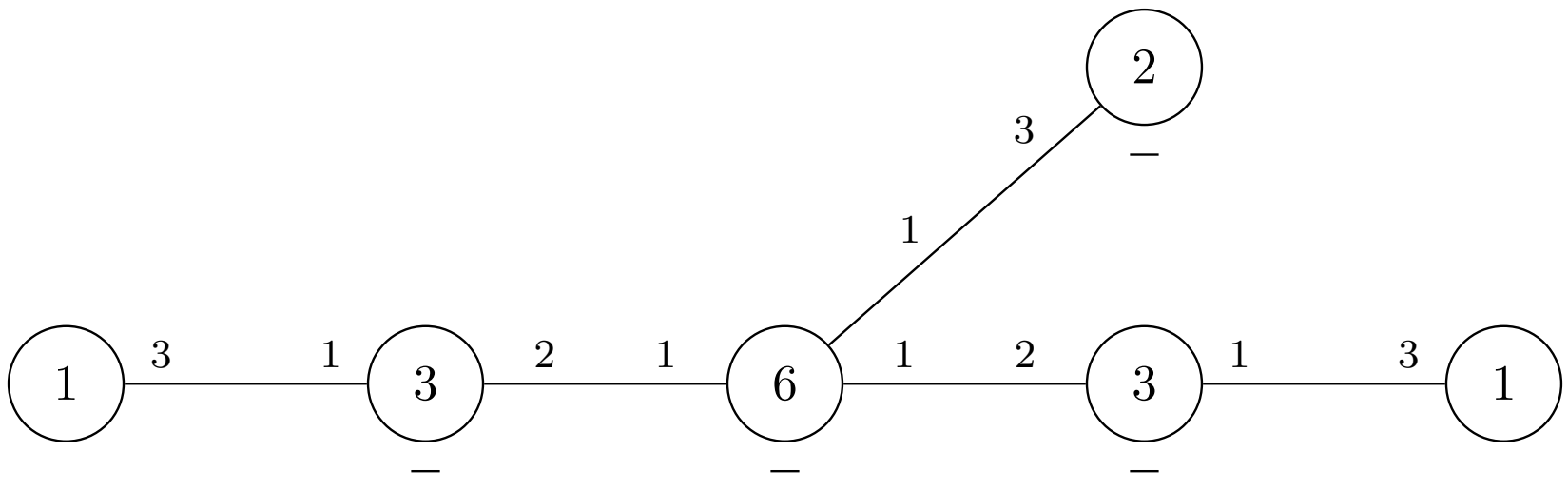
\end{proof}

\begin{proof}[Proof of Theorem~\ref{thm:4k4}]
By Lemma~\ref{lem:v-2k}, $v\leq4k+4$ implies
$k_2\leq2k$. 
If $c_2=1$, then Lemma~\ref{lem:v-c2} applies. 
For $c_2\geq2$, the results follow from Theorem~\ref{thm:2k}. 
\end{proof}

\noindent\textbf{Funding} This work was supported by the National Natural Science Foundation of China (No.~12371338) and the Key Project in Natural
Science Research of Anhui Provincial Department of Education (No.~2023AH050268).\bigskip

\noindent\textbf{Data Availability} No data was used for the research described in the article.

\section*{Declarations}

\textbf{Conflict of interest} The authors declare that they have no known competing financial interests or personal
 relationships that could have appeared to influence the work reported in this paper.


\begin{thebibliography}{99}
\bibitem{BCN89}
A.~E. Brouwer, A.~M. Cohen and A. Neumaier,
\emph{Distance-Regular Graphs}, Springer-Verlag, Berlin, 1989.

\bibitem{BK99}
A.~E. Brouwer and J.~H. Koolen,
The distance-regular graphs of valency four,
\emph{J. Algebraic Combin.} \textbf{10} (1999), 5--24.

\bibitem{CDK13}
M. C\'amara, E.~R. van Dam, J.~H. Koolen and J. Park,
Geometric aspects of $2$-walk-regular graphs,
\emph{Linear Algebra Appl.} \textbf{439} (2013), 2692--2710.

\bibitem{CN07}
M. Conder and R. Nedela,
Symmetric cubic graphs of small girth,
\emph{J. Combin. Theory Ser. B} \textbf{97} (2007), 757--768.

\bibitem{DDF11}
C. Dalf\'o, E.~R. van Dam, M.~A. Fiol, E. Garriga and B.~L. Gorissen,
On almost distance-regular graphs,
\emph{J. Combin. Theory Ser. A} \textbf{118} (2011), 1094--1113. 

\bibitem{DFG09}
C. Dalf\'o, M.~A. Fiol and E. Garriga,
On $k$-walk-regular graphs,
\emph{Electron. J. Combin.} \textbf{16} (2009), 
\#R47. 

\bibitem{DKT16}
E.~R. van Dam, J.~H. Koolen and H. Tanaka,
Distance-regular graphs,
\emph{Electron. J. Combin.} (2016),
\#DS22. 

\bibitem{JKL26}
W. Jin, J.~H. Koolen and C. Lv,
Amply regular graphs with $\mu$ close to half the valency and
group divisible designs, arXiv:2605.25754, 2026.

\bibitem{KP10}
J.~H. Koolen and J. Park,
Shilla distance-regular graphs,
\emph{European J. Combin.} \textbf{31} (2010), 2064--2073.

\bibitem{KP12}
J.~H. Koolen and J. Park,
Distance-regular graphs with $a_1$ or $c_2$ at least
half the valency,
\emph{J. Combin. Theory Ser. A} \textbf{119}
(2012), 546--555.

\bibitem{QDK16}
Z. Qiao, S.~F. Du and J.~H. Koolen,
$2$-Walk-regular dihedrants from group divisible designs,
\emph{Electron. J. Combin.} \textbf{23} (2016), \#P2.51.

\bibitem{QKP16}
Z. Qiao, J.~H. Koolen and J. Park,
$2$-Walk-regular graphs with a small number of vertices compared to the valency,
\emph{Linear Algebra Appl.} \textbf{510} (2016), 10--24.

\bibitem{QPK19}
Z. Qiao, J. Park and J.~H. Koolen,
On $2$-walk-regular graphs with a large intersection number $c_2$,
\emph{European J. Combin.} \textbf{80} (2019), 224--235.

\end{thebibliography}
\end{document}